\documentclass{amsart}
\usepackage{verbatim,amsmath,amssymb,cite,epsfig,xspace,xcolor,mathrsfs,amsthm,booktabs}

\newcommand{\atc}{atomistic-to-continuum }

\newcommand{\C}{\mathcal{C}}
\newcommand{\A}{\mathcal{A}}
\newcommand{\R}{\mathscr{R}}
\newcommand{\lat}{\Lambda}

\newcommand{\Oa}{\Omega^{\rm a}}
\newcommand{\Oc}{\Omega^{\rm c}}

\newcommand{\Ea}{\mathcal{E}^{\rm a}}
\newcommand{\Ec}{\mathcal{E}^{\rm c}}
\newcommand{\Eb}{\mathcal{E}^\beta}

\newcommand{\Y}{\mathscr{Y}}
\newcommand{\U}{\mathscr{U}}
\newcommand{\adm}{\mathscr{A}}

\newcommand{\T}{\mathcal{T}}
\renewcommand{\S}{\mathcal{S}}
\newcommand{\vor}{\mbox{vor}}

\newcommand{\alert}[1]{{\color{black}#1}}
\newcommand{\alerta}[1]{{\color{black}#1}}
\newcommand{\alertco}[1]{{\color{black}#1}}
\newcommand{\alertnew}[1]{{\color{black}#1}}
\newcommand{\ignore}[1]{}

\newcommand{\abs}[1]{|#1|}
\DeclareMathOperator*{\argmin}{argmin}
\newcommand{\Real}{\mathbb{R}}

\newcommand{\dof}{{\rm DoF}}

\newtheorem{remark}{Remark}[section]

\title[Formulation and optimization of the BQCE method]{Formulation and optimization of the
energy-based blended quasicontinuum method}

\author{M. Luskin}
\address{M. Luskin\\ 127 Vincent Hall \\ 206 Church St. SE \\
  Minneapolis \\ MN 55455 \\ USA}
\email{luskin@math.umn.edu}

\author{C. Ortner}
\address{C. Ortner\\ Mathematics Institute \\ Zeeman Building \\
  University of Warwick \\ Coventry CV4 7AL \\ UK}
\email{christoph.ortner@warwick.ac.uk}

\author{B. Van Koten}
\address{B. Van Koten\\ 127 Vincent Hall \\ 206 Church St. SE \\
  Minneapolis \\ MN 55455 \\ USA}
\email{vankoten@math.umn.edu}

\numberwithin{equation}{section}

\date{\today}

\thanks{
ML and BVK were supported in part by the
NSF PIRE Grant OISE-0967140, DOE Award DE-SC0002085, AFOSR Award
FA9550-12-1-0187, and the University of Minnesota Supercomputing
Institute.  CO was supported by EPSRC grant EP/H003096 ``Analysis
  of atomistic-to-continuum coupling methods.''}



\begin{document}
\maketitle

\begin{abstract}
  \alert{We formulate an energy-based atomistic-to-continuum coupling
    method based on blending the quasicontinuum method for the
    simulation of crystal defects. We utilize theoretical results from
    \cite{BvK:blend1d,VKOr:blend2} to derive optimal choices of
    approximation parameters (blending function and finite element
    grid) for microcrack and di-vacancy test problems and confirm our analytical predictions in numerical tests.}
\end{abstract}

\section{Introduction}
A major goal of materials science is to predict the macroscopic
properties of materials from their microscopic structure.  For this
purpose, it is necessary to understand the behavior of defects in
these materials.  We propose a computational tool, the energy-based
blended quasicontinuum method (BQCE), for simulating defects such
as cracks, dislocations, vacancies, and interstitials in crystalline
materials.

Accurate modeling of the region near a defect requires the use of
computationally expensive atomistic models.  Such models are practical
only for small problems.  However, a defect may interact with a large
region of the material through long-range elastic fields.  Thus,
accurate simulation of defects requires the use of a large
computational domain; typically, the size required rules out the use
of atomistic models for the entire region of interest.

Fortunately, the long-range elastic fields generated by a defect are
well described by continuum models which can be efficiently computed
using the finite-element method.  Thus, defects can be accurately and
efficiently simulated by coupled models which use an atomistic model
near the defect and a continuum model elsewhere.  We call any such
model an \emph{\atc coupling}.

Many \atc couplings have been proposed in recent years
\cite{badia:onAtCcouplingbyblending,bauman:applicationofArlequin,
  Ortiz:1995a,Shapeev:2010a,xiao:bridgingdomain, FiKoGuSch:1989,
  IyGa:2011, Shimokawa:2004}; see~\cite{Miller:2008, LiLuOrVK:2011a}
for a survey of \atc couplings and computational benchmark tests.
These couplings fall into two major classes: energy-based and
force-based.  Energy-based couplings provide an approximation to the
atomistic energy of a configuration of atoms.  Force-based couplings
provide a non-conservative force-field which approximates the forces
on each atom under the atomistic model.  Our BQCE method is an
energy-based coupling.  Both types of couplings have intrinsic
advantages; the development of energy-based couplings is especially
important for finite-temperature applications since equilibrium
statistical properties and transition rates can be directly
approximated~\cite{dupuy.finite.qc11,hyperqc11}.

\alert{Here, and throughout, we only consider {\em concurrent}
  coupling methods.  An alternative approach are upscaling methods
  such as ~\cite{bridgingscalepark05}, which achieve
  increased resolution near defects by adding finer scales to a coarse
  scale description rather than by decomposing space into fine
  (atomistic) and coarse (continuum) scale models.}

The primary source of error for most (concurrent) energy-based
couplings is the \emph{ghost force}.  We say that a coupling suffers
from ghost forces if it predicts non-zero forces on the atoms in a
perfect lattice.  Although many attempts have been made to develop an
energy-based coupling free from ghost forces such couplings are
currently known only for a limited range of problems
\cite{Shimokawa:2004, E:2006, Shapeev:2010a, OrtZha:2011a}.

Shapeev's method \cite{Shapeev:2010a} applies to one and
two-dimensional simple crystals with an atomistic energy based on a
pair interaction model, and can be extended to 3D if a modified
``continuum model'' is used \cite{shap3d}.  GR-AC\footnote{\alerta{The
    acronym ``GR-AC'' was introduced in~\cite{OrtZha:2011a}.  It
    stands for ``geometry reconstruction-based atomistic-to-continuum
    coupling method.''}} was proposed in \cite{Shimokawa:2004,
  E:2006}, and has recently been implemented for a two dimensional
crystal with nearest neighbor many-body interactions in
\cite{OrtZha:2011a}.  No ghost force free methods are currently known
for three-dimensional crystals, for multi-lattice crystals (except in
1D \cite{OrSh:multi_pre}), or for atomistic models with general
many-body interactions. The field-based coupling of Iyer and Gavini
\cite{IyGa:2011} is another interesting approach; however, it is
unclear at present whether it is competitive in terms of computational
complexity and generality.

In our BQCE method, the ghost forces cannot be eliminated but can be
controlled in terms of an additional approximation parameter (the
blending width)~\cite{VKOr:blend2, BvK:blend1d}.
BQCE applies to a wide range of problems for which
no ghost force free methods are known; these problems include
three-dimensional crystals with general many-body interactions as
well as multi-lattices.  This makes it an attractive method for such
challenging and physically important problems.

The key feature of BQCE is a \emph{blending region} where the
atomistic and continuum contributions to the total energy are smoothly
mixed.  The ghost forces of the BQCE method can be made arbitrarily
small by increasing the size of this blending region~\cite{VKOr:blend2, BvK:blend1d}.
BQCE shares the idea of a blending region with the bridging domain
method~\cite{xiao:bridgingdomain}, the AtC
coupling~\cite{badia:onAtCcouplingbyblending}, and the
Arlequin~method~\cite{bauman:applicationofArlequin}.  By contrast, the
energy-based quasicontinuum method (QCE) \cite{Ortiz:1995a}, Shapeev's
method \cite{Shapeev:2010a}, and GR-AC \cite{Shimokawa:2004, E:2006,
  OrtZha:2011a} exhibit an abrupt transition between the atomistic and
continuum models.  We call any method with a blending region a
\emph{blended method}, and we call the weights which mix the atomistic
and continuum contributions to the energy a \emph{blending function}.

Both the bridging domain method and AtC coupling are very general
formulations, each of them incorporating BQCE and QCE as special
cases.  Our BQCE formulation provides a set of specific instructions
for the successful implementation of a blended method.  We identify
two important practical differences between BQCE and the bridging
domain and AtC coupling methods.  First, the BQCE method specifies a
strong coupling between the atomistic and continuum regions, whereas
weak couplings based on Lagrange multipliers or the penalty method
have been used in most work involving the bridging domain and AtC
coupling methods.  Second, in BQCE, we blend the atomistic site-energy
with a continuum site-energy based on the continuum site-energy
defined in some formulations of the QCE method (see Section~\ref{sec:
  coarse grained continuum energy}).  This guarantees that BQCE
correctly predicts the total energy of a perfect lattice subjected to
uniform strain.

Our approach to blending is supported by rigorous analysis.  \alert{In
\cite{BvK:blend1d}, we showed that the ghost force error of BQCE in 1D
does indeed decrease with the size of the blending region, and we
found that the error is minimized when the blending function is a Hermite
cubic-spline.}  We also found that the error of the BQCE method in
predicting lattice instabilities can be reduced by increasing the size
of the blending region.

\alertnew{In this paper, we have extended the atom-based formulation
\eqref{eqn: blended quasicontinuum energy} given for 1D problems in
\cite{BvK:blend1d} to a volume-based multiD formulation
\eqref{eq:bqce_practical} that allows finite element coarse-graining
in the continuum region.

The implementation of BQCE requires the choice of two approximation
parameters: a \emph{blending function} $\beta$ and a finite-element
mesh $\T$ which is used to compute the continuum contribution to the
energy. In Section~\ref{sec:complexity}, we give optimal choices of
$\beta$ and $\T$ to minimize global error norms for the problem of a
point defect in a 2D crystal, based on theoretical results in
\cite{VKOr:blend2}.

We demonstrate the validity of these results in computational test
problems in which we simulate a microcrack (see Figure
\ref{fig:mcrack}) and a di-vacancy.}




\subsection*{Outline}
\alert{In Sections~\ref{sec: atomistic energy} and \ref{sec: coarse grained
  continuum energy}, we introduce an atomistic model for simple
lattices with general many-body interactions and its corresponding
continuum Cauchy--Born approximation. In Section~\ref{sec: bqce}, we
give a precise formulation of BQCE in this setting.
In Section~\ref{sec:complexity}, we offer guiding principles on
choosing the blending function $\beta$ and the mesh $\T$ based on our
analysis in~\cite{VKOr:blend2}.  The results of this analysis are
summarized in Table~\ref{tab:complexity bqce}.}
In Section~\ref{sec:num}, we describe the details of our numerical
experiment.  We use an atomistic energy based on the Embedded Atom
Method.  We did not choose our atomistic energy to model any specific
physical material; instead, the atomistic energy is a toy model chosen
for its simplicity.  The observed rates of convergence are in
agreement with the rates predicted in Section~\ref{sec:complexity}.
In particular, the error of BQCE in the $W^{1,2}$ semi-norm decreases
as $\dof^{-\frac{1}{2}}$ where $\dof$ is the number of degrees of
freedom. \alertnew{We summarize our results and discuss open
questions concerning local quantities of interest in the concluding
Section~\ref{conclusion}.}

\section{The Atomistic Energy}\label{sec: atomistic energy}
Let $\lat$ be a $d$-dimensional Bravais lattice.
\alerta{We call $\lat$ the \emph{reference lattice}.
In the present work, we consider only \emph{monatomic crystals}.
That is, we assume that each site of the lattice $\lat$
is the reference position of a single atom, and that all atoms are identical.}
Let
\begin{equation*}
\Y := \{y:\lat\rightarrow\Real^d; y(\xi) \neq y(\eta) \mbox{ for all } \xi \neq \eta \}
\end{equation*}
be the set of \emph{deformations} of $\lat$.  Let $\Oa \subset \lat$
be a finite subset of $\lat$.
We call $\Oa$ the \emph{atomistic computational domain}.

In the Embedded Atom Method (EAM) \cite{Daw:1984a}, the energy of
$\Oa$ subjected to deformation $y$ takes the form
\begin{equation}\label{eqn: embedded atom potential}
\Ea (y) := \sum_{\xi \in \Oa} \bigg \{ \sum_{\substack{\eta \in \lat \\ \eta \neq \xi}} \frac{1}{2}\phi \big(\abs{y(\eta) - y(\xi)}\big)
+ G\bigg(\sum_{\substack{\eta \in \lat \\ \eta \neq \xi}} \rho\big(\abs{y(\eta) - y(\xi)}\big)\bigg) \bigg \},
\end{equation}
where $\phi$ is a pair potential,
$\rho$ is an electron density function,
and $G$ is an embedding function.
We call the inner sum
\begin{equation}
  \label{eq:eam_site_energy}
\Ea_\xi (y) := \sum_{\substack{\eta \in \lat \\ \eta \neq \xi}} \frac{1}{2}\phi \big(\abs{y(\eta) - y(\xi)}\big)
+ G\bigg(\sum_{\substack{\eta \in \lat \\ \eta \neq \xi}} \rho\big(\abs{y(\eta) - y(\xi)}\big)\bigg)
\end{equation}
the \emph{atomistic site-energy of atom $\xi$}.  We observe that the sum
defining the energy of an atom is finite in practice even though
summation ranges over the infinite lattice $\lat$.  This is because
the pair potential $\phi$ and electron density function $\rho$ are
taken to have a \emph{cut-off radius} $r_c$ such that
\begin{equation*}
\phi(r) = \rho(r) = 0 \mbox{ for all } r \geq r_c.
\end{equation*}

For defining the atomistic model and the \alert{BQCE method}, we can
in fact consider a more general class of potentials than EAM.  We only
require that the total energy can be decomposed into a sum of
localized site-energies associated with each atom.  By
\emph{localized}, we mean that the energy associated with an atom
$\xi$ does not depend on the positions of atoms beyond a certain
cut-off distance.  Such an assumption may be violated for certain
energies arising from quantum mechanics, but does hold for most
empirical potentials including EAM, the bond-angle potentials, and so
forth. In addition, we require that the site-energies are {\em
  homogeneous}; that is, the energies of atoms which have the same
local environment are the same.

To make these assumptions precise, we let $\Ea_\xi (y)$ denote the
site-energy associated with atom $\xi$ under deformation $y$ and
require that it is of the form
\begin{equation}\label{eqn: homogeneity of site energy}
\Ea_\xi(y) := V\big( \big\{ y(\eta) - y(\xi) : \eta \in \lat
  \setminus\{\xi\}, |y(\eta) - y(\xi)| < r_c \big\} \big),
\end{equation}
where $V$ is the site-energy potential.
We assume that the resulting
site-energies are twice continuously differentiable, that is, $\Ea_\xi
\in C^2(\Y)$. The restricted dependence of $\Ea_\xi$ on atoms within
the cut-off radius may also be expressed in the form
\begin{equation*}
\frac{\partial \Ea_\xi}{\partial y(\eta)}(y) = 0 \mbox{ for all } \eta \mbox{ with } |y(\eta) - y(\xi)| \geq r_c.
\end{equation*}
This quantifies the requirement that the site-energy is localized.


Given $\Ea_\xi$ we define the energy of $\Oa$ subjected to a
deformation $y$ by
\begin{equation}\label{eqn: general atomistic energy}
\Ea(y) : = \sum_{\xi \in \Oa} \Ea_\xi(y).
\end{equation}
We call an energy of the form~\eqref{eqn: general atomistic energy}
where $\Ea_\xi$ satisfies~\eqref{eqn: homogeneity of site energy}
\emph{homogeneous}.  When we define the Cauchy--Born strain energy
density corresponding to~\eqref{eqn: general atomistic energy} in
Section~\ref{sec: coarse grained continuum energy}, we use that $\Ea$
is homogeneous: if $\Ea$ is not homogeneous, the energy per unit
volume in a perfect lattice subjected to uniform strain \alert{may be
  more difficult to obtain (see \cite{shapML} and references therein
  for examples).}


\begin{remark}
  The locality assumption can, in principle, be replaced by an
  assumption that the interaction strength decays sufficiently rapidly
  with increasing distance between atoms.
\end{remark}
\begin{remark}
  Homogeneity of the site-energy is our main assumption that is
  violated for multi-lattices. We show in \cite{VKOr:blend2} how to
  generalize our formulation for that scenario.
\end{remark}

\section{The Cauchy--Born Site Energy}\label{sec: coarse grained
  continuum energy}

BQCE is a coupling of an atomistic energy based on~\eqref{eqn: general
  atomistic energy} with a coarse-grained continuum elastic energy
based on the Cauchy--Born strain energy density~\eqref{eqn: cauchy
  born strain energy density}.  Let $\vor(\xi)$ denote the Voronoi
cell of a site $\xi \in \lat$.  Then the Cauchy--Born strain energy
density \alert{$W:\Real^{d \times d} \rightarrow \Real \cup \{ \pm
  \infty \}$} corresponding to $V$ is defined by
\begin{equation}\label{eqn: cauchy born strain energy density}
  W(F) := \frac{1}{|\vor(0)|} \Ea_0(y^F),
\end{equation}
where $y^F \in \Y$ is the homogeneous deformation $y^F(\xi) = F\xi$.
$W(F)$ may be interpreted as the energy per unit volume in $\lat$
subjected to strain $F$. Note that the assumption of homogeneity
\eqref{eqn: homogeneity of site energy} ensures that $\Ea_\xi(y^F) =
\Ea_0(y^F)$ for all $\xi \in \lat$.


We will use the Cauchy--Born strain energy density~\eqref{eqn: cauchy
  born strain energy density} to derive a coarse-grained continuum
energy suitable for coupling with the atomistic energy~\eqref{eqn:
  general atomistic energy}.  First, we define a space of
coarse-grained deformations.  Let $\lat^{\rm rep} \subset \lat$ denote
a set of \emph{representative atoms (repatoms)}, let $\T$ be a
triangulation with vertices $\lat^{\rm rep}$, and let $P^1(\T)$ denote
the set of all functions $y:\Real^d \rightarrow \Real^d$ that are
continous and piecewise affine with respect to $\T$.  We call
$P^1(\T)$ the set of \emph{coarse-grained deformations}. The
Cauchy--Born energy of a deformation $y \in P^1(\T)$ in a domain
$\Omega$ is then given by
\begin{displaymath}
  \Ec(y) := \int_{\Omega} W(\nabla y(x)) \,dx.
\end{displaymath}
The Cauchy--Born approximation is analyzed, for example, in
\cite{BLBL:arma2002, E:2007a, Hudson:stab}.

The definition of the QCE method \cite{Ortiz:1995a} and our
construction of the BQCE method in the next section use a Cauchy--Born
site-energy $\Ec_\xi$, which is analogous to the atomistic site-energy $\Ea_\xi$.
For $y\in P^1(\T)$ and $\xi\in \lat$, we define $\Ec_\xi$ by
\begin{align}
\Ec_\xi(y) :&= \int_{\vor(\xi)} W(\nabla y (x)) \, dx \notag \\
&= \sum_{T \in \T} |\vor(\xi) \cap T| W(\nabla y|_T).\label{eqn: continuum energy per atom}
\end{align}
In formula~\eqref{eqn: continuum energy per atom}, $|\vor(\xi) \cap
T|$ denotes the volume of the intersection of $\vor(\xi)$ with the
element $T$.  We observe that the sum on the right hand side of
equation~\eqref{eqn: continuum energy per atom} is finite because only
finitely many elements $T\in\T$ can intersect $\vor (\xi)$.

\section{The Blended Quasicontinuum Energy}\label{sec: bqce}

\subsection{Formulation of the BQCE method}
The BQCE method is an atomistic-to-continuum coupling based on the QCE
method of Tadmor {\it et al.}  \cite{Ortiz:1995a}.  In QCE, the
reference domain $\Oa$ is partitioned into an atomistic region $\A$
and a continuum region $\C$, and the QC energy $\mathcal{E}^{\rm
  qc}:P^1(\T) \rightarrow \Real$ is defined by
\begin{equation}\label{eqn: quasicontinuum energy}
\mathcal{E}^{\rm qc}(y) := \sum_{\xi \in \A} \Ea_\xi(y) + \sum_{\xi \in \C} \Ec_\xi(y). 
\end{equation}

In BQCE, the atomistic and continuum energies per atom are weighted
averages. Given a {\em blending function} $\beta: \Oa \rightarrow
[0,1]$ the BQCE energy $\Eb: P^1(\T) \rightarrow \Real$ is defined by
\begin{equation}\label{eqn: blended quasicontinuum energy}
\Eb(y) := \sum_{\xi\in\Oa} \beta(\xi) \Ec_\xi(y)+(1 - \beta(\xi)) \Ea_\xi(y).
\end{equation}
We observe that the QCE energy with continuum region $\C$ is the same
as the BQCE energy with $\beta$ chosen as the characteristic function
of $\C$.  Our formulation of BQCE is similar in spirit to the bridging
domain method~\cite{xiao:bridgingdomain}, the AtC
coupling~\cite{badia:onAtCcouplingbyblending}, and the
Arlequin~method~\cite{bauman:applicationofArlequin}; \alert{it differs
  from \cite{xiao:bridgingdomain, badia:onAtCcouplingbyblending,
    bauman:applicationofArlequin} in that we specify a strong coupling
  between the atomistic and continuum models, and accordingly our
  method is based on the minimization of an energy. Moreover, our
  formulation in terms of atomistic and continuum energies per atom
  guarantees that the BQCE energy is exact under homogenous
  deformations.}

The BQCE energy can be rewritten in the form
\begin{align}
  \notag
\Eb(y) &= \sum_{\xi \in \Oa}  \beta(\xi) \int_{\vor(\xi)} W(\nabla y)
\, dx  + (1-\beta(\xi))\Ea_\xi(y) \\
\notag
&= \sum_{\xi \in \Oa} \sum_{T \in \T} \beta(\xi) |\vor(\xi) \cap T|
W(\nabla y|_T) +  \sum_{\xi \in \Oa}(1-\beta(\xi))\Ea_\xi(y)\\
\label{eq:bqce_practical}
&= \sum_{T \in \T} v^\beta_T W(\nabla y|_T) + \sum_{\xi \in \Oa} (1-\beta(\xi))\Ea_\xi(y),
\end{align}
where the \emph{BQCE-effective volume} of the element $T$ is defined by
\begin{equation}
v^\beta_T := \sum_{\xi \in \Oa} \beta(\xi)|\vor(\xi) \cap T|.
\end{equation}

\ignore{
Let $f:\Oa \rightarrow \Real^d$ be an external force.
We define $\Eb_{total} : P^1(\T) \rightarrow \Real$ by
\begin{equation*}
\Eb_{total} (y) := \Eb(y) + \sum_{\xi \in \Oa} f(\xi) \cdot y(\xi).
\end{equation*}
The minimization problem for the BQCE energy is to find
\begin{equation}
 y \in \argmin_{z \in \adm} \Eb_{total}(z),
\end{equation}
where $\adm = \adm(E, F, \T)$ for periodic boundary conditions
or $\adm = \adm(y_0, \Omega^{rep}, \T)$ for Dirichlet boundary conditions.
}

\begin{remark}
  \label{rem: triangulation need not cover lattice}
The triangulation $\T$ need not cover the entire domain $\Oa$.
For $\T$ a triangulation which covers only part of $\Real^d$,
define
\begin{align*}
\Omega(\T) :&= \cup_{T \in \T} T, \text{ and } \\
P^1(\T) :&= \{y: \Omega(\T) \cup \lat \rightarrow \Real^d: \text{ y piecewise affine w.r.t. $\T$ on $\Omega(\T)$}\}.
\end{align*}
We observe that $\Eb(y)$ is defined for $y \in P^1(\T)$ if for every $\xi \in \Oa$ such that $\beta(\xi) > 0$ we have $\emph{\vor} (\xi) \subset \Omega(\T)$.
In particular, it is not necessary to assume that the triangulation $\T$ is refined to atomistic scale anywhere in the domain $\Oa$.
It is possible that the use of a mesh which is not refined to atomistic scale may make the implementation of BQCE easier and more efficient in some cases.
\end{remark}

\alert{
\begin{remark}[Multi-lattices]
  If one interprets a multi-lattice as a simple lattice with a basis
  (see, e.g., \cite{DoElLuTa:2007, VKOr:mle, VKOr:blend2}), then the
  formulation of the atomistic energy \eqref{eqn: general atomistic
    energy} and of the BQCE energy \eqref{eqn: blended quasicontinuum
    energy} require no changes. The only difference is that in the
  general multi-lattice case, the site energy $\Ea_\xi(y)$ has fewer
  symmetries (however, see\cite{VKOr:mle} for examples with high
  symmetry).
\end{remark}
}

\subsection{Far-field boundary conditions}\label{subsec: far field bcs}
A typical application of the BQCE method is the simulation of a defect
or defect region in an infinite crystal. To that end, we require
far-field boundary conditions at the domain boundary. We propose two
choices.

\subsubsection{Dirichlet boundary conditions}
Let $\S \subset \T$ be a finite subset of $\T$, and let $\Omega(\S) :=
\cup_{S \in \S} S$ be a polygonal domain.  When Dirichlet boundary
conditions are imposed, the deformation of the boundary $\partial
\Omega(\S)$ of $\Omega(\S)$ is fixed to agree with some $y_0 \in
P^1(\T)$.  Precisely, we let
\begin{equation}
{\rm Adm} := \big\{y \in P^1(\T): y(x) = y_0(x) \mbox{  for all } x \in \partial \Omega(\S)\big\}
\end{equation}
denote the space of admissible deformations. We then solve the problem
\begin{equation}\label{eqn: minimization problem}
\mbox{Find } y \in \argmin_{z \in {\rm Adm}} \Eb(z),
\end{equation}
where we interpret $\argmin_{z \in {\rm Adm}} \Eb(z)$ as the set of
\emph{local} minimizers of $\Eb$.

\subsubsection{Periodic boundary conditions}
A popular method to construct artificial far-field boundary conditions
is to formulate the problem in a periodic cell. To that end, suppose
that $\Omega(\T) = \Real^d$, $\T$ is periodic, and let $\S\subset\T$ be
the finite element mesh on one periodic cell.
That is, suppose that for some nonsingular matrix $A \in \Real^{d \times d}$
\alerta{whose columns are elements of $\lat$}, we have
\begin{displaymath}
  \T = \bigcup_{n \in \mathbb{Z}^d} \big \{ A n + \S \big\},
\end{displaymath}
and that the union is disjoint.

Let $F \in \Real^{d \times d}$ be nonsingular,
and call $F$ a \emph{macroscopic strain}.
Given $F$, we define the admissible set as
\begin{equation}
  \label{eq:per_bc}
  {\rm Adm} := \big\{y \in P^1(\T): y(x + A n) = y(x) + F A n \mbox{  for
    all } x \in \Real^d \big\}.
\end{equation}
The associated variational problem can again be stated as \eqref{eqn:  minimization problem}.

\alerta{The columns of the matrix $A$ should be understood as the
  directions in which the displacement $u(x) = y(x) - F x$ is
  periodic.  Equivalently, the set $\big \{Ax : x \in [0,1]^d \big \}$
  is a periodic cell.}  \alerta{The macroscopic strain $F$ should be
  understood as the average strain imposed on the crystal.  If the
  size of the periodic cell is large, periodic boundary conditions
  with macroscopic strain $F$ approximate far-field boundary
  conditions with uniform strain $F$ imposed at infinity.}

\alert{Note that, while we only discuss a static formulation, we see no
  obstacle that prevents both our formulation and analysis to apply
  also to {\em quasi-static} problems, in which the energy is
  minimized at each loading step. However, we stress that further
  careful analysis is required to understand the accuracy, in
  particular stability, of the BQCE method as the deformation
  approaches bifurcation points (e.g., as in the formation or motion
  of crystal defects or the propagation of a crack).}

\alertco{
\subsection{Degrees of freedom (DoF)}
In Section \ref{sec:complexity}, we present error estimates for the
BQCE method in terms of the {\em computational complexity}, by which
we mean the computational cost to compute the BQCE energy and its
gradient. This is proportional to the number of degrees of freedom of
the space of admissible functions ${\rm Adm}$. Similarly, in Section
\ref{sec:num}, we plot errors against the number of degrees of
freedom in the computed solution.

The number of degrees of freedom (DoF) is simply the size of the
nonlinear system that we solve to minimize the BQCE energy. In terms
of the notation introduced in the previous section, this is the number
of finite element nodes in the triangulation $\mathcal{S}$ plus the
number of lattice sites in $\lat \setminus \Omega(\mathcal{S})$ for
which $\beta(\xi) < 1$.

In the reduced atomistic computation, which we describe in Section
\ref{sec:num:atm}, DoF denotes the number of unconstrained lattice
sites.  }

\subsection{Computational complexity and optimal parameters}
\label{sec:complexity}
In~\cite{VKOr:blend2}, we conjecture an error estimate for BQCE in 2D
and 3D.  The conjecture is based on our error analysis of BQCE in
1D~\cite{BvK:blend1d} and on the consistency estimates for BQCE in 2D
and 3D in~\cite{VKOr:blend2}.  Following
\cite[Sec. 7.1]{OrtShap:2011a}, we use our conjectured error estimate
to derive \emph{complexity estimates}, which are bounds on the error
of the BQCE method in terms of the number of degrees of freedom.  We
use our complexity estimates to guide the choice of optimal
approximation parameters $\T$ and $\beta$.

Let $y_{\rm a}$ be a \alerta{local minimizer} of the atomistic energy,
and let $y_\beta$ be a \alerta{local minimizer} of BQCE with the same
boundary conditions as $y_{\rm a}$.  Let $h$ be the mesh size function
of the triangulation $\T$ ($h(x) = {\rm diam}(T)$ for a.e. $x \in T$),
and let $\beta$ be the blending function.  Then we conjecture an
estimate of the form:
\begin{equation}\label{eqn: conjectured error estimate}
\begin{split}
{\rm Err} := \|\nabla y_{\rm a} - \nabla y_\beta \|_{L^p} &\lesssim \|h \nabla^2 y_{\rm a} \|_{L^p(C)} + \|\nabla^2 \beta \|_{L^p} \\
    &=: {\rm CG} + {\rm GF},
\end{split}
\end{equation}
\alerta{for given $p \in [1,\infty]$.
For $p=2$, we have proven this result in~\cite{VKOr:blend2};
for $p \neq 2$, it would be difficult to establish this result rigorously;
and for $p \in \{1,\infty\}$, it is unclear in what generality the result holds.}
In~\eqref{eqn: conjectured error estimate}, $\nabla^2 y_{\rm a}$ and $\nabla^2
\beta$ should be interpreted as the second derivatives of smooth
interpolants of $y_{\rm a}$ and $\beta$, and $C := {\rm supp}(\beta)$.
The first term, ${\rm CG} = \|h \nabla^2 y_{\rm a} \|_{L^p(C)}$, is
the finite element {\em coarsening error}, while the second term,
${\rm GF} = \|\nabla^2 \beta \|_{L^p}$, measures the effect of the ghost
forces.

We now describe the problem of a point defect in a 2D crystal.
To quantify the notion of a point defect, we assume that for some $\alpha > 0$,
\begin{equation}
|\nabla^2 y_{\rm a}(x)| \simeq r^{-\alpha}, \qquad \text{where } r = |x|.
\end{equation}
It has been observed in numerical experiments that $\alpha = 2$ for a
dislocation, and that $\alpha = 3$ for a vacancy \cite{FrankMerwe1949,
  OrtShap:2011a}.

We assume that the reference domain $\Omega$ is a
roughly circular region of radius $N$ atomic spacings centered at the
origin.  We let $K_0 > 0$ be the radius of the atomistic region
surrounding the defect, and we let $K_1 > 0$ be the width of the
blending region.  We then choose a radial blending function $\beta$ of
the form
\begin{equation}\label{eq: form beta for complexity}
\beta(x) :=
\begin{cases}
0  &\mbox{if }  r < K_0, \\
\beta_0 \left (\frac{|x| - K_0}{K_1} \right ) &\mbox{if } K_0 \leq r < K_0+K_1,\\
1 &\mbox{if } K_0+K_1 \leq r,
\end{cases}
\end{equation}
where $\beta_0: [0,1] \rightarrow [0,1]$ is a twice continuously differentiable function with $\beta_0(0) = \beta_0'(0) = \beta_0'(1)= 0$ and $\beta_0(1) = 1$.

We summarize our complexity estimates for BQCE in Table~\ref{tab:complexity bqce}.
These estimates are proved in~\cite{VKOr:blend2}.
We distinguish three cases based on the value of $\gamma := \frac{\alpha p}{p + 2}$.
In the second column, we give the optimal rates of convergence for BQCE.
The optimal rates are attained
when $K_1$ is given in terms of $K_0$ by the formula appearing in the third column
and when the mesh size function $h$ is given by
\begin{equation}\label{eqn: mesh size function}
h(x) = \left ( \frac{|x|}{K_0}\right )^\gamma.
\end{equation}
\alerta{
\begin{remark}
When $\gamma > 1$ and $\alpha >2$,
another choice of parameters results in convergence at the optimal rate
but with a better constant: we use
\begin{equation*}
K_1 = K_0^\mu \mbox{\quad for \quad} \mu = \frac{\alpha -\frac{2}{p}}{2 - \frac{2}{p}}
\end{equation*}
in our numerical experiments below.
This choice guarantees that the finite element error
decreases at the same rate as the ghost force error.
See~\cite{VKOr:blend2} for additional suggestions related to choosing parameters.
\end{remark}
}

\alerta{
\begin{remark}
To make the explicit calculation of optimal parameters possible,
we have made simplifying assumptions on the forms of $\beta$ and $h$.
In~\eqref{eq: form beta for complexity}, we assume that $\beta$ is radial.
We explain how to construct $\beta$ given general atomistic, continuum, and blending regions
in Section~\ref{subsec: setup of bqce} below.
Similarly, we assume that $h$ is continuous
even though a nonconstant continuous function cannot be the mesh size function of a triangulation.
In practice, we implement BQCE using a triangulation whose mesh size function approximates
the continuous $h$ defined in~\eqref{eqn: mesh size function};
we explain how to do this in Section~\ref{subsec: setup of bqce}.
We expect that the complexity estimates in Table~\ref{tab:complexity bqce}
will hold for practical choices of $\beta$ and $h$,
not only for the simplified $\beta$ and $h$ used in the derivation of the estimates.
\end{remark}}

\begin{table}\label{table1}
\caption{Complexity estimates and optimal approximation parameters for BQCE~\cite{VKOr:blend2}.
The estimates above are for the problem of a point defect in a two-dimensional crystal.
The variables $K_0$, $K_1$, $N$, and $\alpha$ are defined in Section~\ref{sec:complexity},
and $\gamma := \frac{\alpha p}{p + 2}$.
In all cases, the optimal rate of convergence is attained when $h =\left ( \frac{|x|}{K_0}\right )^\gamma$
and when $K_1$ is given in terms of $K_0$ by the formula in the third column above.
}
\begin{tabular}{c@{\qquad}l@{\qquad}l@{\qquad}}
\toprule
Case & Complexity Estimate & Optimal Parameters \\
\midrule
$\gamma > 1$ & $\| \nabla y_{\rm a} - \nabla y_\beta \|_{L^p} \lesssim \dof^{\max\{ \frac{1}{p} - 1, \frac{1}{p} - \frac{\alpha}{2}\}}$ & $K_1 = K_0$ \\
\midrule
$\gamma = 1$ & $\| \nabla y_{\rm a} - \nabla y_\beta \|_{L^p} \lesssim \dof^{\max\{\frac{1}{p}-1, -\frac{1}{2}\}}$ &
$K_1 = K_0 \ln \left ( \frac{N}{K_0} \right )^{\frac{1}{2}}$\\
\midrule
$\gamma < 1$& $\| \nabla y_{\rm a} - \nabla y_\beta \|_{L^p} \lesssim \dof^{\max\{\frac{1}{p}-1, -\frac{1}{2}\}}$ &
$K_1 = K_0^\gamma N^{1 - \gamma}$\\
\bottomrule
\end{tabular}
\label{tab:complexity bqce}
\end{table}

\ignore{
\begin{table}\label{tab:complexity shapeev}
\begin{tabular}{c@{\qquad}l@{\qquad}l}
\toprule
Case & Complexity Estimate \\
\midrule
$\gamma > 1$ & $\| \nabla y_{\rm a} - \nabla y_\beta \|_{L^p} \lesssim \dof^{\frac{1}{p} - \frac{\alpha}{2}}$ \\
\midrule
$\gamma = 1$ & $\| \nabla y_{\rm a} - \nabla y_\beta \|_{L^p} \lesssim \dof^{-\frac{1}{2}}$ \\
\midrule
$\gamma < 1$& $\| \nabla y_{\rm a} - \nabla y_\beta \|_{L^p} \lesssim \dof^{-\frac{1}{2}}$  \\
\bottomrule
\end{tabular}
\caption{Complexity estimates for Shapeev's method.
These estimates appear in Section~7.1 of~\cite{OrtShap:2011a}.
As in Table~\ref{tab:complexity bqce}, $\gamma := \frac{\alpha p}{p + 2}$. }
\end{table}}

\begin{remark}
  The rates of convergence depend on the geometry of the problem and
  on the norm in which the error is measured.
  We observe that the error of BQCE does not decrease with $\dof$ when measured in the
  $W^{1,1}$-seminorm; \alerta{see the case $p=1$ in Table~\ref{tab:complexity bqce}.}
  This is because the $W^{-1,1}$-norm of the
  ghost force does not decrease as the size of the blending region
  increases.
  \alerta{We refer the reader to~\cite{VKOr:blend2} for further discussion of the
  convergence of BQCE for other geometries and in other norms.}
\end{remark}

\begin{remark} \alertnew{We have proved that linear blending functions have
a reduced rate of convergence compared to smooth blending functions
of the form~\eqref{eq: form beta for complexity} (see~\cite{BvK:blend1d} for 1D
and \cite{VKOr:blend2} for 2D).  In Table~\ref{tab:complexity lin}, we give
comparative convergence rates for
the problem of a point defect in
a two-dimensional crystal
with decay $\alpha=3.$ We emphasize, however, that this is only an
asymptotic estimate. Indeed in our micro-crack experiment in
Section~\ref{sec:num:rates}, we observe a significant preasymptotic regime
where the two blending functions yield comparable errors.
}
\begin{table}
\begin{tabular}{c|c|c}
\toprule
 & linear BQCE & smooth BQCE \\
 & $\| \nabla y_{\rm a} - \nabla y_\beta \|_{L^p}$ & $\| \nabla y_{\rm a} - \nabla y_\beta \|_{L^p}$
 \\
\midrule
$p=2$ &  $\lesssim
\dof^{-\frac14}$ & $\lesssim
\dof^{-\frac12}$  \\
\midrule
$p=\infty$ &  $\lesssim
\dof^{-\frac12}$ & $\lesssim
\dof^{-1}$  \\
\bottomrule
\end{tabular}
\vspace{.2in}

\caption{Complexity estimates comparing the optimized linear
and smooth BQCE method
for the problem of a point defect in
a two-dimensional crystal
with
decay $\alpha=3.$
}
\label{tab:complexity lin}
\end{table}
\end{remark}

\begin{remark}
Although our subsequent numerical
  experiments are restricted to defects with zero Burgers vector
  \alertnew{and consequently $\gamma>1,$} we
  consider more general situations in the above analysis, in order
  to emphasize the generality of our approach. The implementation of
  a/c couplings for defects with non-zero Burgers vector poses
  additional challenges (what should one use for the reference
  domain in the neighborhood of a single dislocation?) that we will address in future work. As a
  matter of fact, we are not aware of numerical results from a/c
  coupling methods, simulating a single dislocation core in an
  otherwise defect-free lattice (as opposed to the simulation of
  a dislocation pair with total zero Burgers vector).
  \end{remark}
\subsection{Comparison of BQCE with other methods}
We use the estimates above to compare the complexity of BQCE
with direct atomistic simulation and with patch test consistent couplings,
such as Shapeev's method.
\subsubsection{Shapeev's method and other patch test consistent methods}
  When $\alpha$ is small (slow decay of the elastic field), BQCE
  converges at the same rate as patch test consistent methods such as
  Shapeev's method.  In fact, if $\alpha \leq 2$ (e.g., a
  dislocation), then the $W^{1,2}$-error of Shapeev's method decreases
  as $\dof^{-\frac{1}{2}}$ with the number of degrees of
  freedom~\cite{OrtShap:2011a}.  We predict the same rate of
  convergence for BQCE when $\alpha \leq 2$.  On the other hand, when
  $\alpha$ is larger, patch test consistent couplings may converge
  faster than BQCE.  When $\alpha > 2$ (e.g., a vacancy, microcrack,
  or dislocation dipole), the $W^{1,2}$-error of Shapeev's method
  decreases as $\dof^{\frac{1}{2} - \frac{\alpha}{2}}$, whereas the
  $W^{1,2}$-error of BQCE decreases as $\dof^{-\frac{1}{2}}$.  Roughly
  speaking, this is due to the fact that, for small $\alpha$, the
  coarse-graining error dominates, whereas for large $\alpha$, the
  ghost force error dominates.

  According to the estimates above, Shapeev's method always converges at least as fast as BQCE.
  However, we remind the reader that patch test consistent methods are currently known
  only for a limited range of problems~\cite{Shimokawa:2004, E:2006, Shapeev:2010a, OrtZha:2011a}.
  As far as we are aware, BQCE is the only convergent energy-based method for 2D or 3D crystals
  with general many-body interactions.
  Moreover, even when patch test consistent methods are known,
  BQCE may be an attractive method for slowly decaying defects,
  since it converges at the same rate as patch test consistent methods for these problems.
  See~\cite{VKOr:blend2} for a more detailed analytic comparison of BQCE with Shapeev's method.

\alert{\subsubsection{Reduced atomistic computations (ATM)}
\label{sec:num:atm}
In our numerical experiments below, we compare the performance of BQCE
with the \emph{reduced atomistic computation (ATM)}, in which
far-field boundary conditions are simulated by imposing Dirichlet
boundary conditions directly on the atomistic model.  Precisely, let
$\Omega \subset \lat$ be an \emph{atomistic computational domain}, and
let $F$ be a macroscopic strain as in Section~\ref{subsec: far field
  bcs}.  We then define the admissible space
\begin{equation*}
{\rm Adm} := \big\{y: \lat \rightarrow \Real^d: y(\xi) = F\xi \mbox{  for all } \xi \in \lat \setminus \Omega \},
\end{equation*}
and we solve the variational problem
\begin{equation}\label{eqn: atomistic minimization problem}
\mbox{Find } y \in \argmin_{z \in {\rm Adm}} \Ea(z).
\end{equation}

We regard ATM as the simplest means of simulating far-field boundary
conditions, and we feel that an atomistic-to-continuum coupling is
successful only when it is more accurate than ATM.
In~\cite{VKOr:blend2}, we show that
\begin{equation*}
\big \|\nabla y^\Omega_{\rm a} - \nabla y_{\rm a} \big \|_{L^p} \lesssim N^{\frac{2}{p} + 1 -\alpha} \simeq \dof^{\frac{1}{p} + \frac{1}{2}- \frac{\alpha}{2}},
\end{equation*}
where $N$ is the diameter of $\Omega$, $y^\Omega_{\rm a}$
solves~\eqref{eqn: atomistic minimization problem}, ${p \in
  [1,\infty]}$, and $\alpha > 2$.  Comparing the estimate above with
Table~\ref{tab:complexity bqce}, we observe that BQCE converges more
quickly than ATM when $\alpha < 3$, at the same rate when $\alpha =
3$, and more slowly when $\alpha > 3,$ which we summarize in
Table~\ref{tab:complexity atm}. Our numerical experiments
concern the case $\alpha = 3$, and we do observe that BQCE and ATM
converge at the same rate.  However, for optimally chosen parameters,
BQCE is still much more accurate than ATM by a significant constant
factor.
}

\alertnew{
\begin{table}
\begin{tabular}{c|c|c|c|c}
\toprule
 & BQCE & ATM &$2<\alpha$ &ATM rate vs. \\
 & $\| \nabla y_{\rm a} - \nabla y_\beta \|_{L^p}$ & $\| \nabla y_{\rm a} - \nabla y_{\rm a}^\Omega \|_{L^p}$
 &&BQCE rate\\
\midrule
&  &  &$\alpha>3$&higher\\
$p=2$&  $\lesssim
\dof^{-\frac12}$ & $\lesssim
\dof^{-\frac12+\frac12(3-\alpha)}$&$\alpha=3$ &same \\
&  & &$\alpha<3$&lower \\
\midrule
&  &  &$\alpha>3$&higher\\
$p=\infty$ &  $\lesssim
\dof^{-1}$ & $\lesssim
\dof^{-1+\frac12(3-\alpha)}$&$\alpha=3$&same \\
&  & &$\alpha<3$&lower \\
\bottomrule
\end{tabular}
\vspace{.2in}

\caption{Complexity estimates for the optimized BQCE
and ATM
for the problem of a point defect in
a two-dimensional crystal
with
decay $\alpha>2$~\cite{VKOr:blend2}.
The variable $\alpha$ is defined in Section~\ref{sec:complexity}.
}
\label{tab:complexity atm}
\end{table}
}

\section{Numerical Example}
\label{sec:num}
\subsection{Setup of the atomistic model}
\begin{figure}
  \begin{center}
    \includegraphics[height=4.5cm]{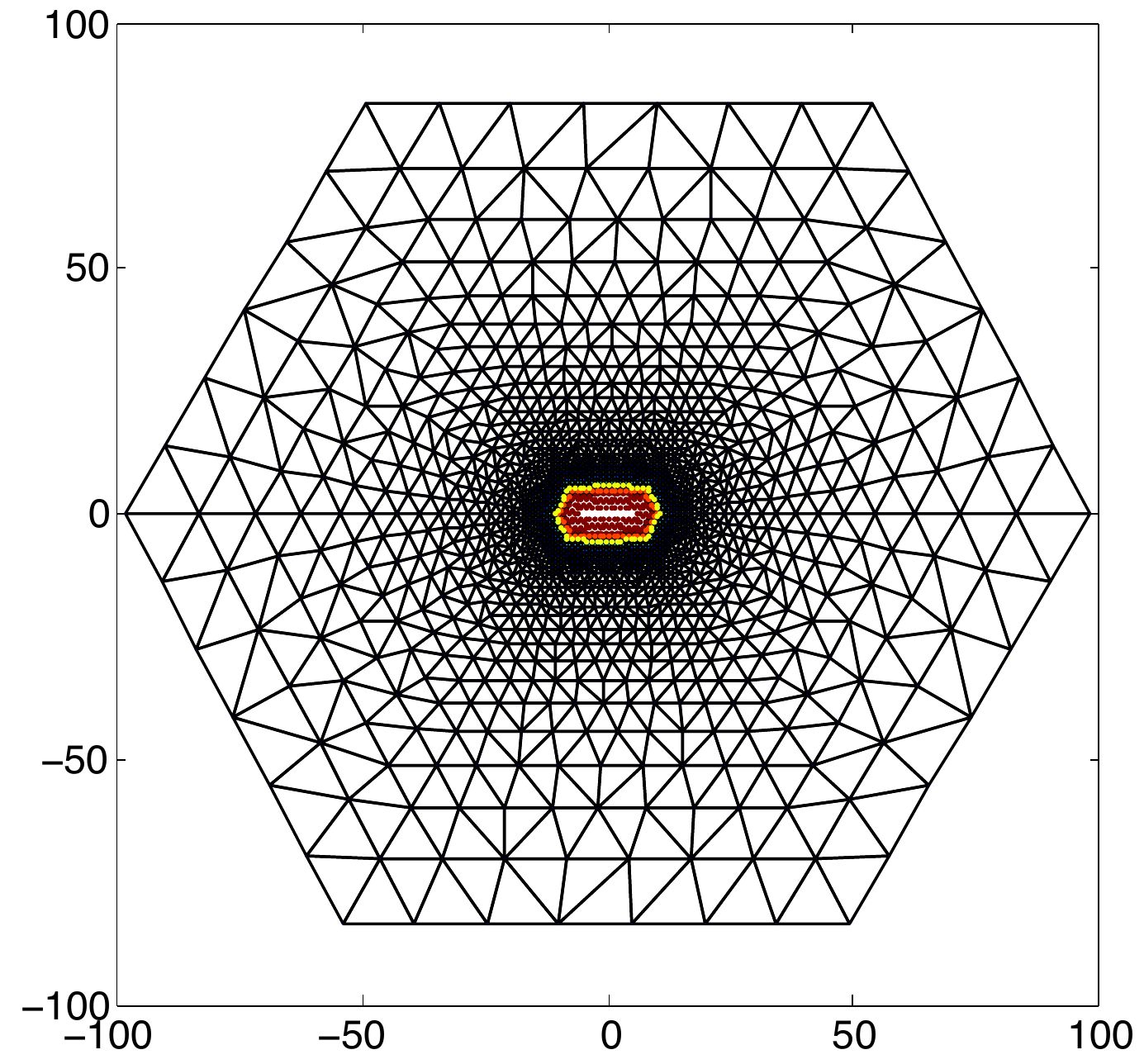}
    \qquad
    \includegraphics[height=4.5cm]{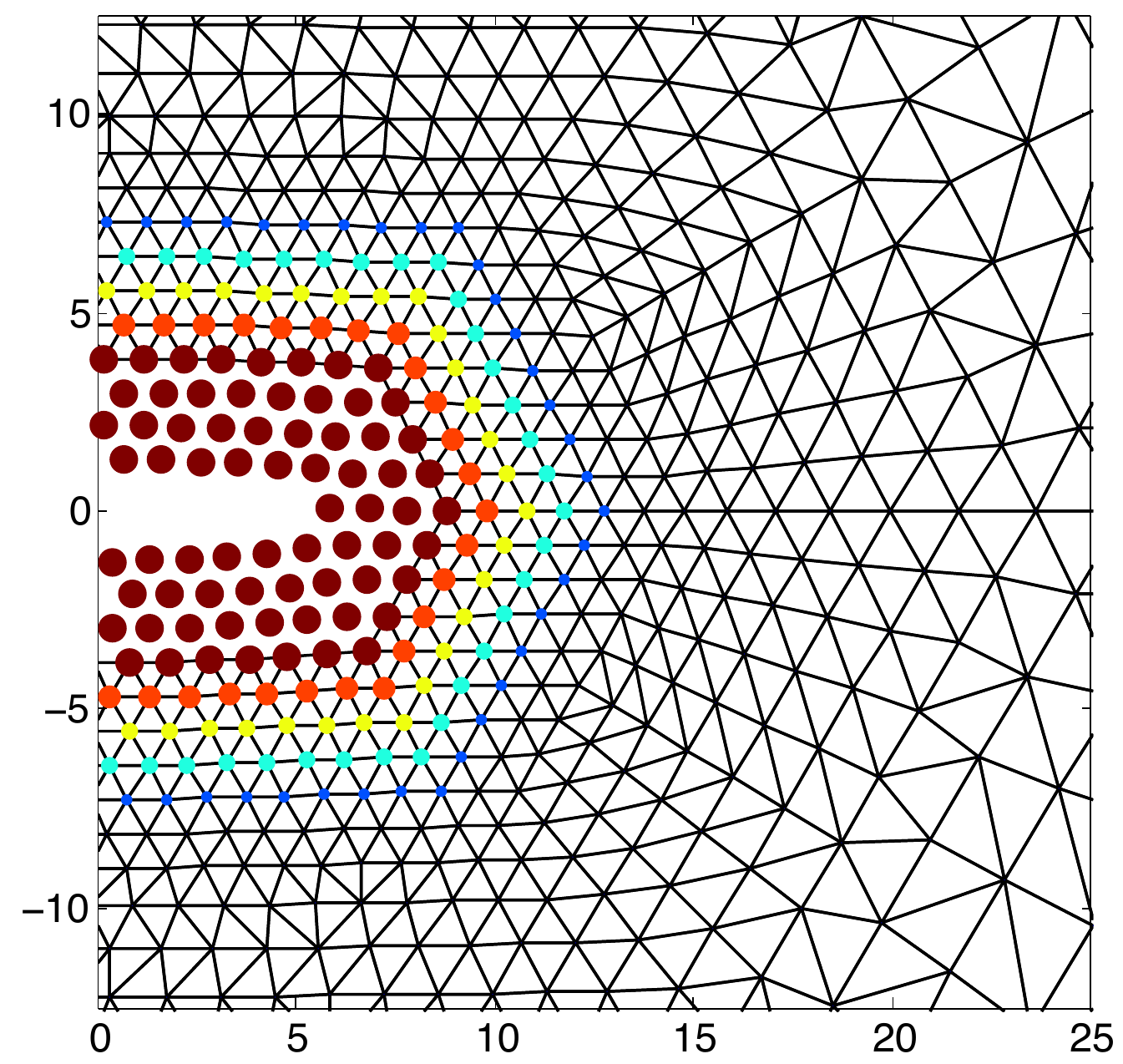}
  \end{center}
  \caption{\label{fig:mcrack} Deformed configuration in atomic units
    of a BQCE solution for a microcrack in a computational domain with
    approximately $3 N^2$ atoms ($N = 100$ in this figure, but $N =
    500$ in the benchmarks described in Sections \ref{sec:num:rates}
    and \ref{sec:num:divac}). The color and size of the atom positions
    indicate the value of the blending function.}
\end{figure}

In the 2D triangular lattice defined by
\begin{displaymath}
  \lat :=\left(\begin{array}{cc} 1 &1/2 \\ 0 &
      \sqrt{3}/2 \end{array}\right) \cdot \mathbb{Z}^2
\end{displaymath}
we choose a hexagonal domain $\Oa$ with embedded microcrack as
described in Figure~\ref{fig:mcrack}. The sidelength of the domain is
\alertco{$N = 500$} atomic spacings, and a defect is introduced by
removing a segment of atoms at the center of $\Oa$,
\alertco{$\lat^{\rm def}_{11} := \{ -5 e_1, \dots, 5 e_1 \}$ in the
  microcrack example in Section \ref{sec:num:rates}, and $\lat^{\rm
    def}_2 := \{ 0, e_1 \}$ in the di-vacancy example in Section
  \ref{sec:num:divac},} \alertnew{where the subscripts refer to the number
  of atoms removed in the simulation.}

\alertco{We supply the domain with Dirichlet boundary conditions by setting
$y_{\rm a}(\xi) = F \xi$ for $\xi \in \lat \setminus \Oa$, where $F$
is a macroscopic strain that we will specify below.  }


The site energy is given by the EAM toy-model
\eqref{eq:eam_site_energy}, with
\begin{align*}
  \phi(r) =~& \big[ e^{-2 a (r-1)} - 2 e^{- a (r-1)} \big], \quad
  \rho(r) = e^{-b r}, \quad \text{and}  \\
  G(\bar\rho) =~& c \big[ (\bar\rho-\bar\rho_0)^2 + (\bar\rho -
  \bar\rho_0)^4 \big],
\end{align*}
where $a, b, c, \bar\rho_0 \in \mathbb{R}$ are parameters of the
model.  In all our computational experiment we choose
\begin{displaymath}
  a = 4.4, \quad b = 3, \quad c = 5, \quad \bar\rho_0 = 6 e^{-b}, \quad
  r^{\rm cut}_1 = 1.8, \quad r^{\rm cut}_2 = 2.5.
\end{displaymath}
\alertco{The model mimics qualitative properties of practical EAM
  models such as those described in \cite{eam_johnson, eam_mishin}.

To simplify the implementation of the model, we have
  introduced a cut-off in the reference configuration instead of the
  deformed configuration. In \eqref{eq:eam_site_energy}, we compute
  only sums over first, second and third neighbors. Since there are
  no rearrangements of atom connectivity in our experiments, we do not
  expect that this affects the outcome of the results. }


\subsection{Setup of the BQCE method}\label{subsec: setup of bqce}
Our implementation of the BQCE method is based on
\eqref{eq:bqce_practical}, using standard finite element assembly
techniques. The construction of the blending function is governed by
two approximation parameters:
\begin{itemize}
\item $K_0 \in \mathbb{N}$ denotes the number of atomic layers
  surrounding the microcrack where $\beta = 1$;
\item $K_1 \in \mathbb{N}$ denotes the number of atomic layers in the
  blending region.
\end{itemize}

For the microcrack problem we expect $\alpha = 3$ in the context of
Section \ref{sec:complexity}. Hence, according to Table
\ref{tab:complexity bqce}, the choice $K_1 := K_0$ (so that the number
of atoms in the atomistic region and blending region are comparable)
results in the optimal rate of convergence for both $p = 2$ and $p = \infty$.

Let $d_{\rm hop}(\xi, \eta)$ denote the hopping distance in the
triangular lattice (with natural extension to sets), then we define
\begin{align*}
  \lat^{\rm a} :=~& \big\{ \xi \in \Oa : d_{\rm hop}(\xi, \lat^{\rm
    crack}) \leq K_0 \big\}, \quad \text{and} \\
  \lat^{\rm b} :=~& \big\{ \xi \in \Oa : K_0 < d_{\rm hop}(\xi,
  \lat^{\rm crack}) \leq K_0 + K_1 \big\}.
\end{align*}
We consider three choices of the blending function, which are all
easily defined for general interface geometries:
\begin{itemize}
\item {\it QCE:} choosing $\beta$ to be the characteristic function of
  the atomistic region $\lat^{\rm a}$, and $K_1 = 0$, yields the QCE
  method defined in \eqref{eqn: quasicontinuum energy}.
\item {\it Linear Blending:} Let $d(\xi)$ denote the hopping distance from the
  atomistic region, then we choose
  \begin{displaymath}
    \beta_{\rm lin}(\xi) := \max(1, d(\xi) / K_1).
  \end{displaymath}
\item {\it Smooth Blending:} \alerta{The three nearest-neighbour
    lattice vectors are given by
\begin{displaymath}
a_1 := \left(\begin{array}{c} 1 \\ 0 \end{array}\right), \quad
a_2 := \left(\begin{array}{c} 1/2 \\ \sqrt{3}/2 \end{array}\right),
\mbox{ and } a_3 := \left(\begin{array}{c} -1/2 \\ \sqrt{3}/2 \end{array}\right).
\end{displaymath}}
Define
$\Delta_i^2 \beta(\xi) := \beta(\xi+a_i) -
  2 \beta(\xi) + \beta(\xi-a_i)$, and
\begin{displaymath}
\Phi(\beta) := \sum_{\xi
    \in \Oa} \sum_{i = 1}^3 |\Delta_i^2 \beta(\xi)|^2.
\end{displaymath}
Then we define
  \begin{equation*}
    \qquad\quad\beta_{\rm smooth} := {\rm argmin} \big\{ \Phi(\beta) : \beta(\xi)
    = 0 \text{ in $\lat^{\rm a}$ and } \beta(\xi) = 1 \text{ in $\Oa
      \setminus \lat^{\rm a} \cup \lat^{\rm b}$} \big\}.
  \end{equation*}
  Roughly speaking, this choice generalizes the cubic spline blending
  suggested in \cite{BvK:blend1d}. It provides a practical variant for
  interface geometries that are not circular.
\end{itemize}

The third approximation parameter is the finite element mesh in the
continuum region. We coarsen the finite element mesh away from the
boundary of the blending region according to the rule suggested by the
complexity estimates in Table \ref{tab:complexity bqce}. As a matter
of fact it turns out that the mesh size growth is too rapid to create
shape-regular meshes, hence we also impose the restriction that
neighboring element can at most grow by a prescribed factor; this
introduces an additional logarithmic factor in the complexity
estimates \cite[Sec. 7.1]{OrtShap:2011a}.

The resulting energy functional is minimized using a preconditioned
Pol\'{a}k--Ribi\`{e}re conjugate gradient algorithm described in
\cite{LiLuOrVK:2011a}. We removed the termination criterion in this
algorithm and allowed it to converge to its maximal precision, that
is, until the numerically computed descent direction ceases to be an
actual descent direction for the energy; this occurs at an accuracy of around
$10^{-5}$ in atomic units. \alertco{We then perform three Newton
  steps, after which the residual is in the range of machine
  precision, to ensure that the results are not affected by the
  accuracy of the nonlinear solver.}

\subsection{Microcrack}
\label{sec:num:rates}
\alertco{In the microcrack experiment, we remove a long segment of atoms,
$\lat^{\rm def}_{11} = \{-5 e_1, \dots, 5 e_1\}$ from the
computational domain. The body is then loaded in mixed mode I \& I\!I,
by setting
\begin{displaymath}
  F := \left(\begin{array}{cc} 1 & \gamma_{\rm I\!I} \\ 0 & 1 +
      \gamma_{\rm I} \end{array}\right) \cdot F_0,
\end{displaymath}
where $F_0 \propto I$ minimizes the Cauchy--Born stored energy
function $W$, $\gamma_{\rm I} = \gamma_{\rm I\!I} = 0.03$ ($3\%$ shear
and $3\%$ tensile stretch).}

We solve the BQCE problem (to be precise, the reduced atomistic
problem (ATM) for increasing domain sizes, QCE and BQCE problems for
linear and smooth blending functions) for increasing parameters $K_0$,
and compute the error relative to the exact atomistic solution.

\begin{figure}
  \begin{center}
    \includegraphics[width=9cm]{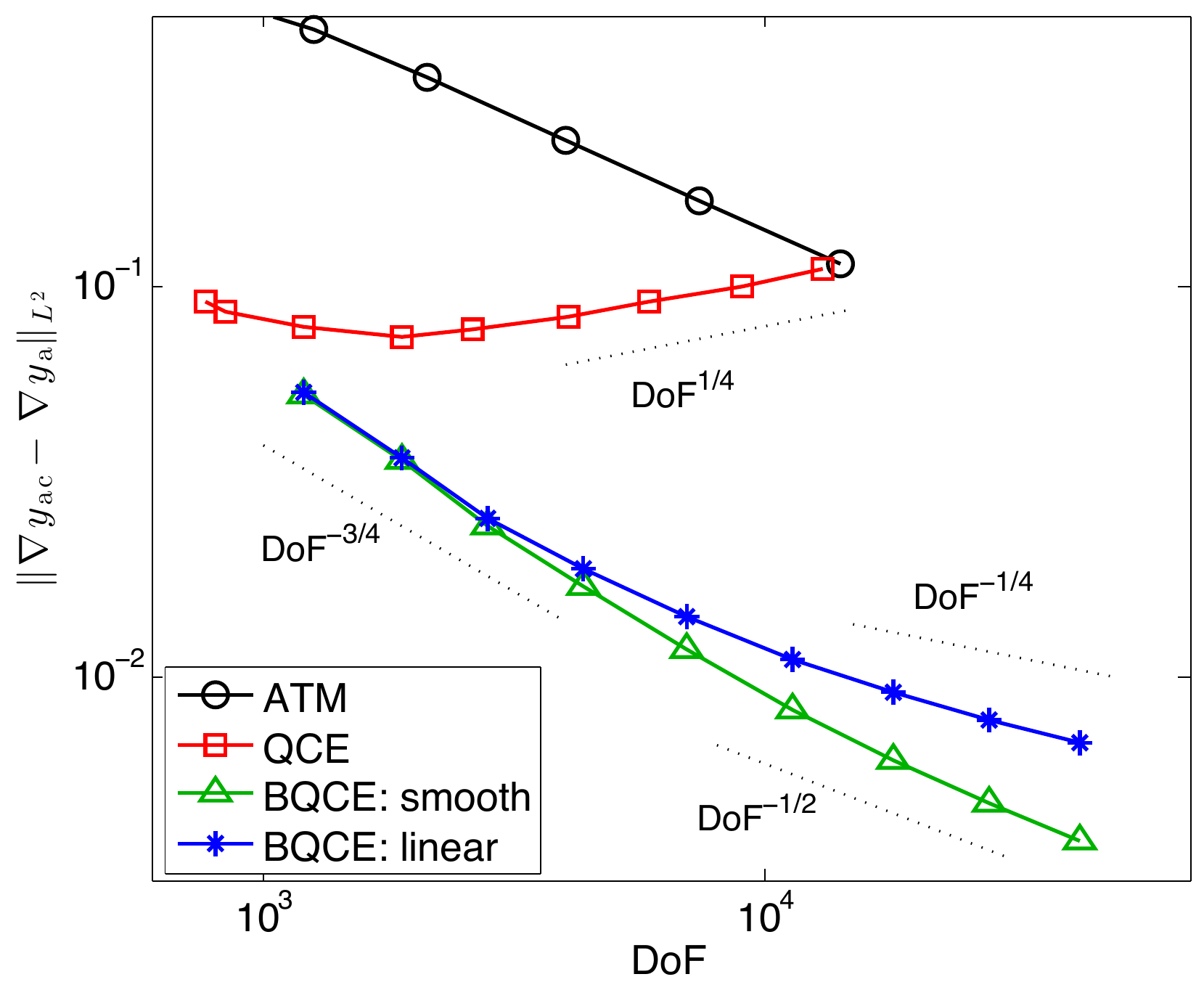}
  \end{center}
  \caption{\label{fig:mcrack_err2} Convergence rates in the
    energy-norm (the $H^1$-seminorm) for the microcrack benchmark
    problem described in Section \ref{sec:num:rates}.  This corresponds to
    the case $p = 2$ and $\alpha = 3$ in Table~\ref{tab:complexity
      bqce}.  In the above, $y_{\rm a}$ denotes the minimizer of the
    atomistic energy, and $y_{\rm ac}$ denotes the corresponding
    minimizer of the ATM, QCE or BQCE problems.}
\end{figure}

\begin{figure}
  \begin{center}
    \includegraphics[width=9cm]{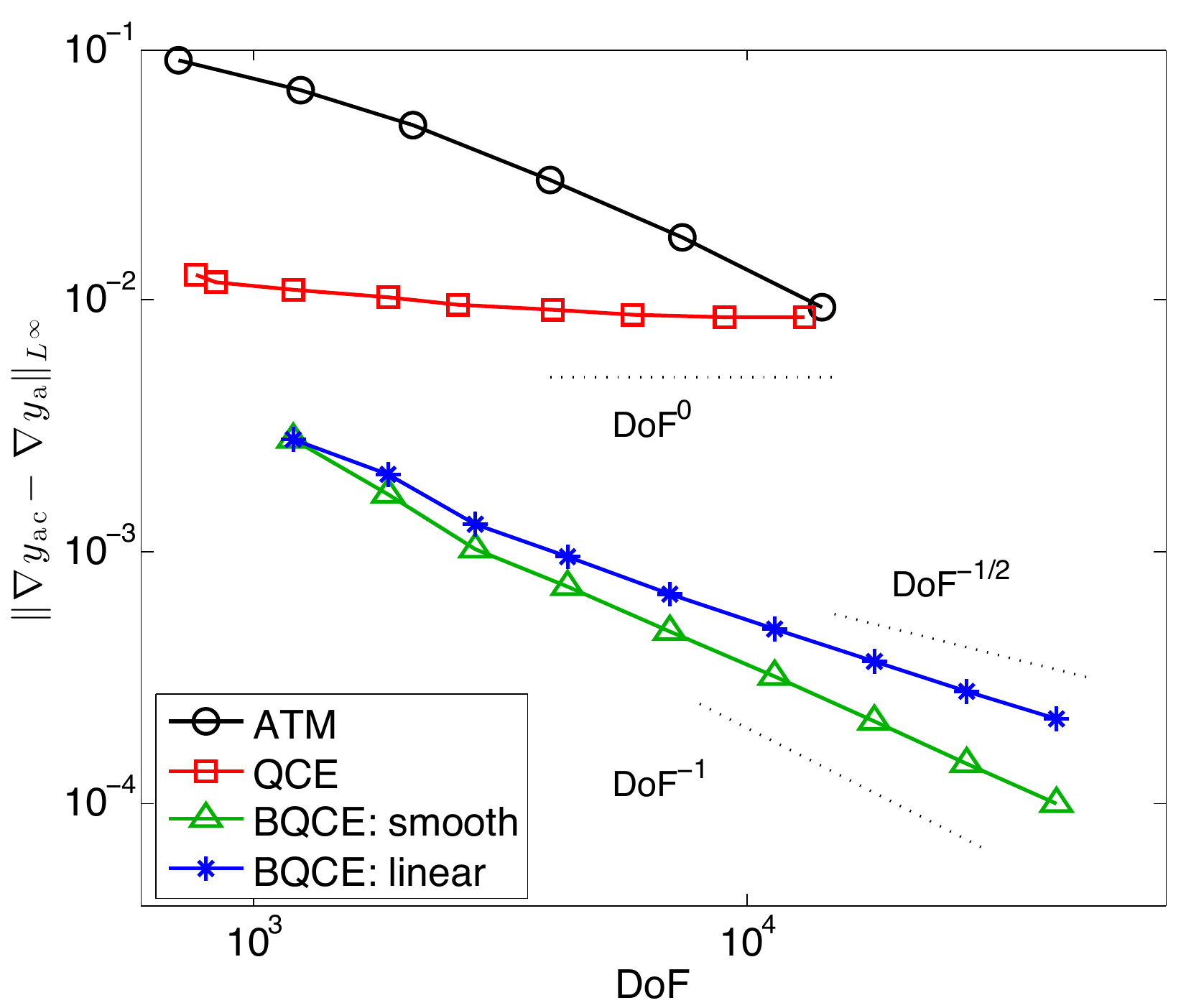}
  \end{center}
  \caption{\label{fig:mcrack_errinf} Convergence rates in the
    $W^{1,\infty}$-seminorm for the microcrack benchmark problem
    described in Section \ref{sec:num:rates}.  This corresponds to the case
    $p = \infty$ and $\alpha = 3$ in Table~\ref{tab:complexity bqce}.
    In the above, $y_{\rm a}$ denotes the minimizer of the atomistic
    energy, and $y_{\rm ac}$ denotes the corresponding minimizer of
    the ATM, QCE or BQCE problems.}
\end{figure}

\begin{figure}
  \begin{center}
    \includegraphics[width=9cm]{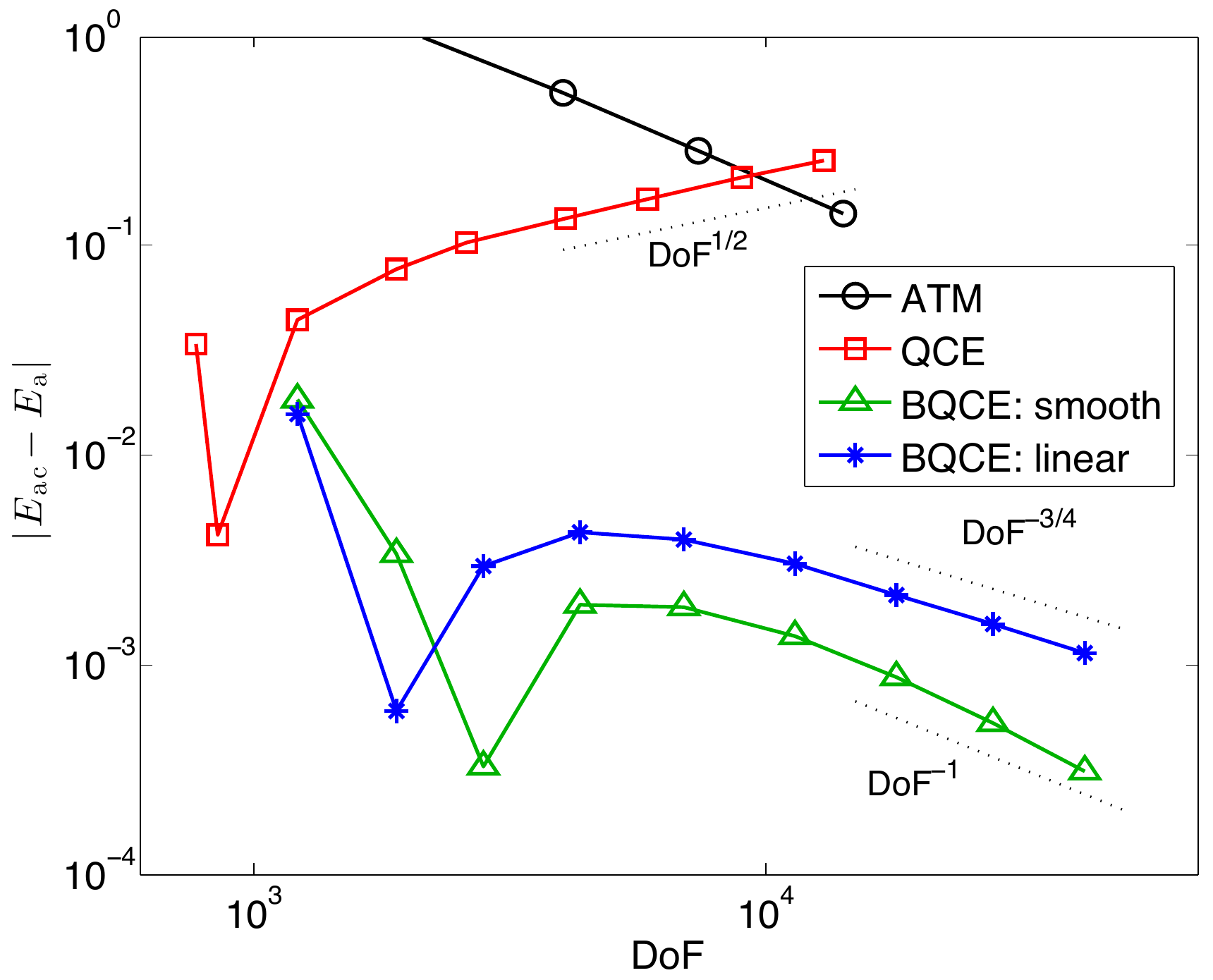}
  \end{center}
  \caption{\label{fig:mcrack_errE} Convergence rates in the energy,
    for the microcrack benchmark problem described in Section
    \ref{sec:num:rates}.  In the above, $E_{\rm a}$ denotes the minimum of
    the atomistic energy, and $E_{\rm ac}$ denotes the minimum of
    either the ATM, QCE or BQCE energies.}
\end{figure}

The relative errors in the $W^{1,2}$-seminorm are displayed in Figure
\ref{fig:mcrack_err2}; the relative errors in the
$W^{1,\infty}$-seminorm are displayed in Figure
\ref{fig:mcrack_errinf}; the errors in the energy are displayed in
Figure \ref{fig:mcrack_errE}. \alertco{We observe a significant pre-asymptotic
regime where the rate of convergence is faster than predicted in our
theory. For larger computations, the rate of convergence approaches
closely the predicted rate. }

In particular, it is worth noting that the advantage of a smooth
blending function only becomes significant in the asymptotic regime
and is less pronounced than our theory might suggest. Since the
precomputation of $\beta_{\rm smooth}$ is computationally cheap and
straightforward, we nevertheless propose the choice $\beta_{\rm
  smooth}$ for most simulations.

\alertco{Finally, we remark that the initial dip in
the energy error in Figure~\ref{fig:mcrack_errE} is most likely due to a change of sign in the energy
error. }

\alertco{
\subsection{Di-vacancy}
\label{sec:num:divac}
We perform a second numerical experiment, in order to demonstrate that
the effects we observed in the microcrack example, which are not
covered by our theory, are indeed caused by a pre-asymptotic
regime. To this end we consider a di-vacancy defect, where only two
neighboring sites $\lat^{\rm def}_2 := \{0, e_1\}$ are removed from
the lattice. We apply 3\% isotropic stretch and 3\% shear loading, by
setting
\begin{displaymath}
  F := \left(\begin{array}{cc} 1 + s & \gamma_{\rm I\!I} \\ 0 & 1 +
      s \end{array}\right) \cdot F_0,
\end{displaymath}
where $F_0$ minimizes $W$, $s = \gamma_{\rm I\!I} = 0.03$.

The relative errors in the $W^{1,2}$-seminorm are displayed in Figure
\ref{fig:divac_err2}; the relative errors in the
$W^{1,\infty}$-seminorm are displayed in Figure
\ref{fig:divac_errinf}; the errors in the energy are displayed in
Figure \ref{fig:divac_errE}. \alertnew{For the $W^{1,2}$ and
$W^{1,\infty}$-seminorms we now observe precisely the rates of
convergence ${\rm DoF}^{-1/2}$ and ${\rm DoF}^{-1}$
predicted by our theory. However, we observe a convergence rate ${\rm DoF}^{-3/2}$
for the energy, rather than the square of the rate for the
$W^{1,2}$-seminorm which is ${\rm DoF}^{-1}$.}
We can,
at present offer no rigorous explanation for this improved convergence
rate, but speculate that it is caused by a cancellation effect that is
not captured by our theory.

\begin{figure}
  \begin{center}
    \includegraphics[width=9cm]{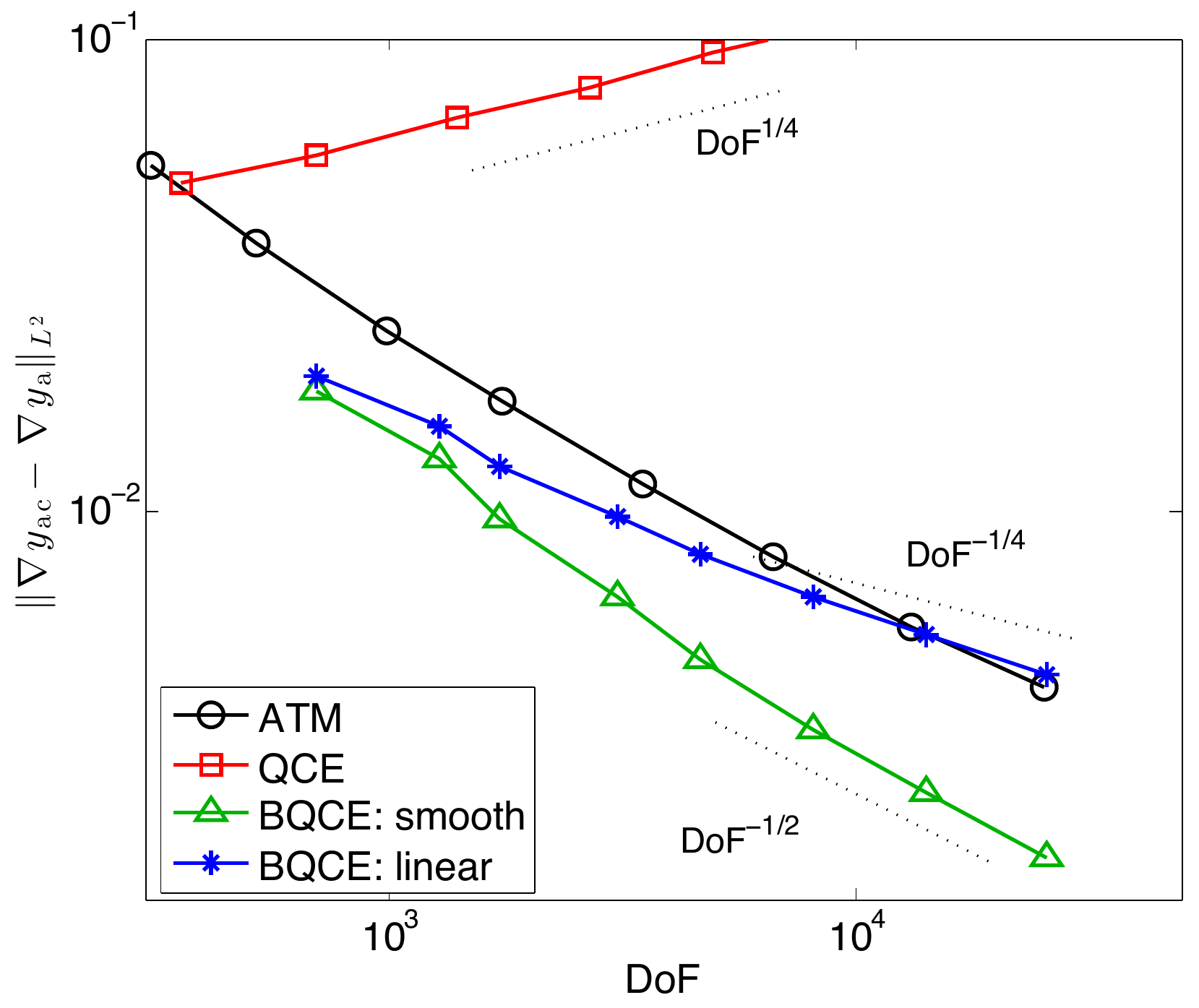}
  \end{center}
  \caption{\label{fig:divac_err2} Convergence rates in the energy-norm
    (the $H^1$-seminorm) for the di-vacancy benchmark problem
    described in Section \ref{sec:num:divac}.  This corresponds to the
    case $p = 2$ and $\alpha = 3$ in Table~\ref{tab:complexity bqce}.
    In the above, $y_{\rm a}$ denotes the minimizer of the atomistic
    energy, and $y_{\rm ac}$ denotes the corresponding minimizer of
    either the ATM, BQCE or QCE problems.}
\end{figure}

\begin{figure}
  \begin{center}
    \includegraphics[width=9cm]{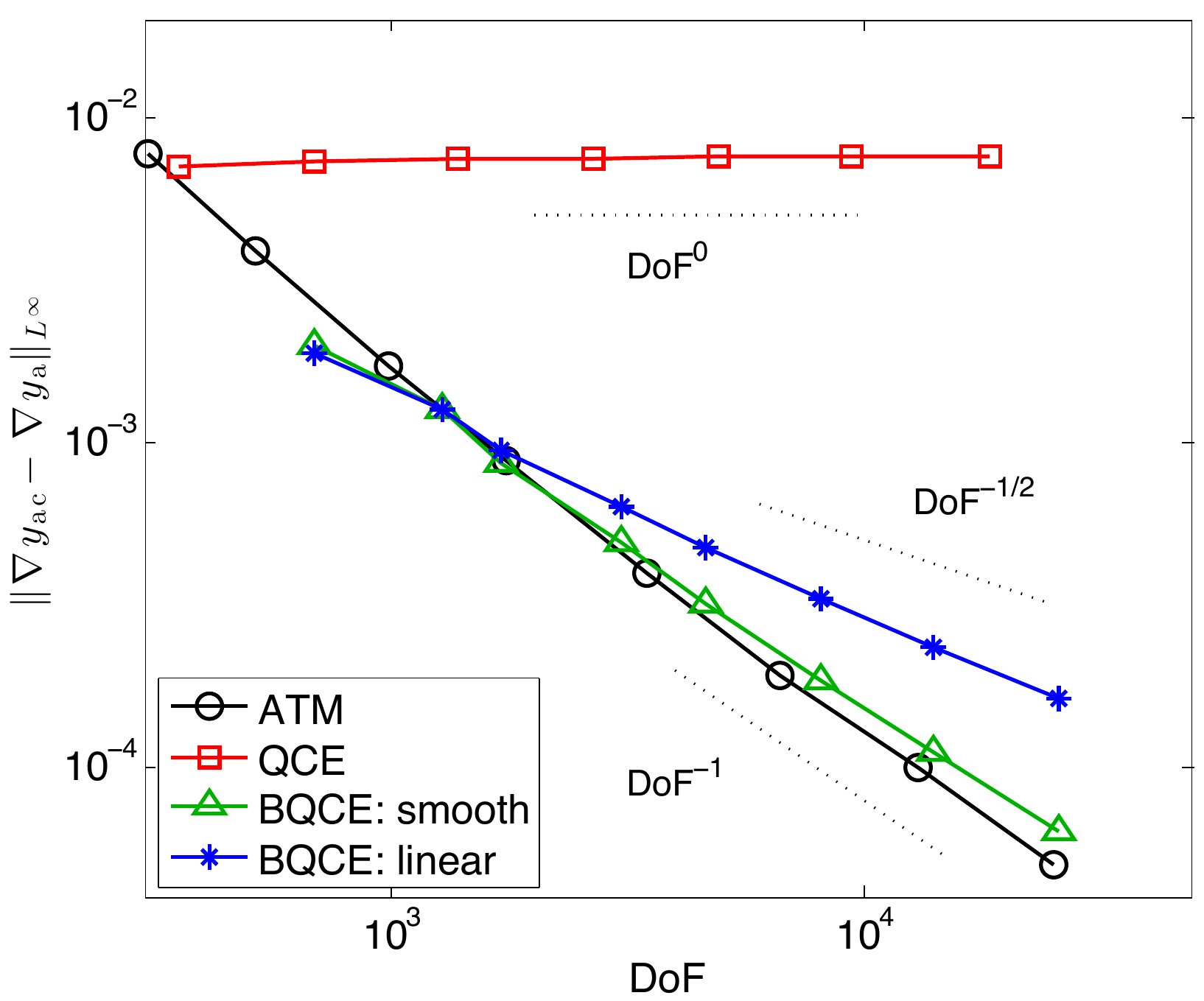}
  \end{center}
  \caption{\label{fig:divac_errinf} Convergence rates in the
    $W^{1,\infty}$-seminorm for the di-vacancy benchmark problem
    described in Section \ref{sec:num:divac}.  This corresponds to the
    case $p = \infty$ and $\alpha = 3$ in Table~\ref{tab:complexity
      bqce}.  In the above, $y_{\rm a}$ denotes the minimizer of the
    atomistic energy, and $y_{\rm ac}$ denotes the corresponding
    minimizer of the ATM, BQCE or QCE problems.}
\end{figure}

\begin{figure}
  \begin{center}
    \includegraphics[width=9cm]{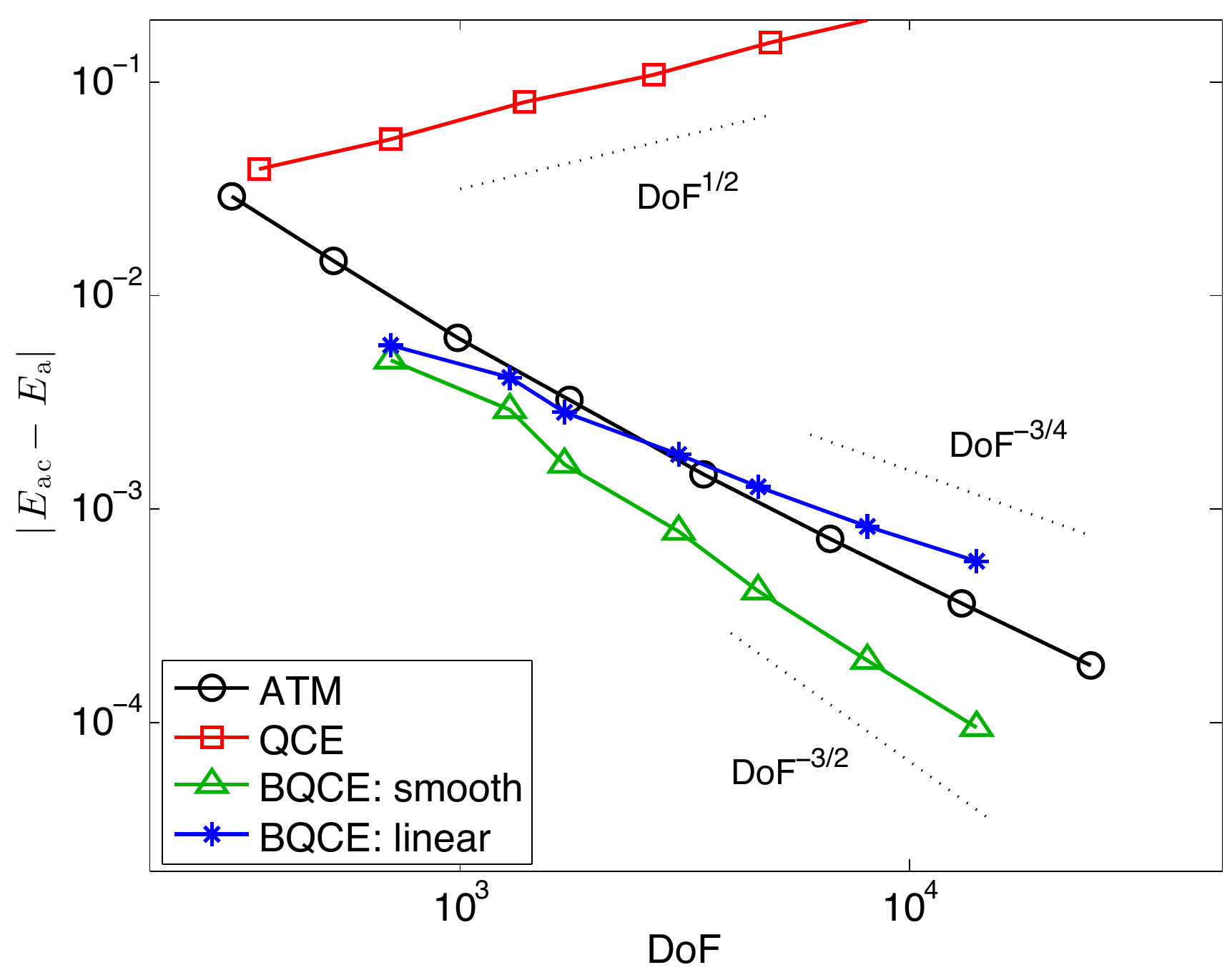}
  \end{center}
  \caption{\label{fig:divac_errE} Convergence rates in the energy,
    for the di-vacancy benchmark problem described in Section
    \ref{sec:num:divac}.  In the above, $E_{\rm a}$ denotes the minimum of
    the atomistic energy, and $E_{\rm ac}$ denotes the corresponding
    minimum of either the ATM, QCE or BQCE energies.}
\end{figure}

It is important to point out that, in this second numerical
experiment, BQCE is {\em not} substantially more efficient than {\rm
  ATM}. In fact, in the $W^{1,\infty}$-seminorm the errors are almost
identical. This example clearly demonstrates that an atomistic
coarse-graining scheme need not always lead to improved computational
results, but that the improvement over a judiciously performed
atomistic simulation is highly problem dependent. We also remark that
for similar examples (e.g., a single vacancy with $F = F_0$), we
observed that ATM has a higher accuracy than BQCE in both the
$W^{1,2}$ and $W^{1,\infty}$-seminorms.  }

\section{Conclusion}
\label{conclusion}

\alert{We have formulated an atomistic-to-continuum coupling, which we
  call the \break
  blended energy-based quasicontinuum method (BQCE).  Our
  formulation requires only two approximation parameters for the
  implementation of BQCE: the blending function $\beta$ and the finite
  element mesh $\T$. For the problems of a microcrack and a di-vacancy in a
  two-dimensional crystal, we utilized theoretical results from
  \cite{VKOr:blend2} to obtain optimal choices of approximation
  parameters (blending function and finite element grid) to minimize a
  global error norm and confirmed our analytical predictions in
  numerical tests.

  An interesting open question is how well different
  atomistic-to-continuum couplings (in particular, QCE, BQCE and the
  patch test consistent methods) perform if the quantity of interest
  is localized in the atomistic region.
  In~\cite{PrBDBaElOd:2008}, this question was investigated through a series of numerical
  experiments in which the Arlequin method was used to simulate a one-dimensional chain with a defect.
  The authors concluded that the local error of the Arlequin method
  could be controlled by moving the blending region away from the defect
  and increasing its size.  For the quasicontinuum method, local quantities
  of interest were investigated in \cite{arndtluskin07b,arndtluskin07c,OdenPrudhommeBauman:2006}.
  We hope to develop rigorous error estimates and optimal approximation parameters
  for localized quantities of interest in a future work.



  At present BQCE is a competitive choice for numerical simulations of
  atomistic multi-scale problems as it applies to a much broader
  class than any known patch test consistent method.  For example,
  Shapeev's method applies only to two-dimensional, simple lattice
  crystals with pair interactions.  Our BQCE method can be applied to
  challenging and physically important problems featuring
  three-dimensional, multi-lattice crystals and arbitrary many-body
  interaction potentials.   }

\end{document}